\documentclass[a4paper,10pt]{article}

\def \Z {{\mathbf Z}}
\def \R {{\mathbf {R}}}

\def\uu{\bigsqcup}
\def\eps{\varepsilon}
\title{ \bf Around Kalikow's theorem on multiple mixing}
\author{\bf Valery V. Ryzhikov}
\date{}
\usepackage{graphicx}
\graphicspath{{pictures/}}
\DeclareGraphicsExtensions{.pdf,.png,.jpg}
\begin{document}

\maketitle

\vspace{20mm} 
\Large

Some results in ergodic theory on joinings and multiple mixing 

and the proof of Steve Kalikow's theorem on 3-fold mixing 

of 2-fold mixing rank one automorphisms   are discussed.

\newpage

\begin{center}
\includegraphics{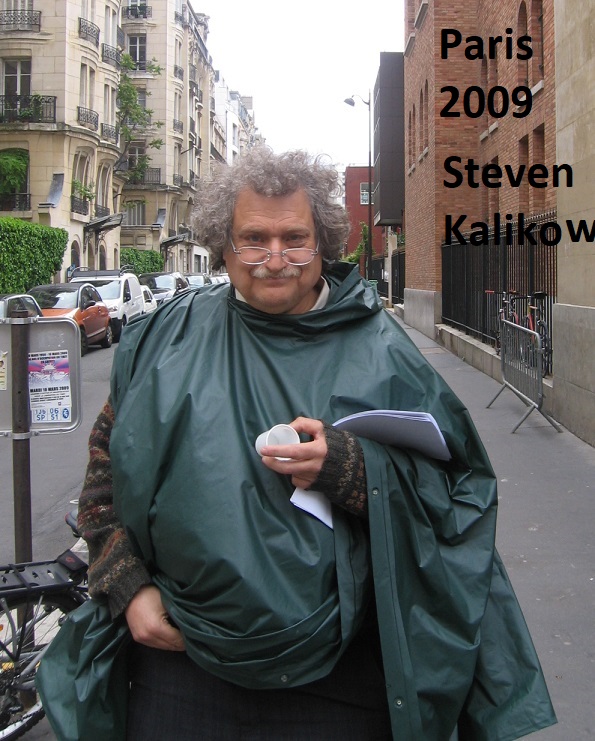}
\end{center}

\Large

Kalikow's theorem \cite{Ka} on 3-mixing of 2-mixing automorphisms of rank 1
was preceded by Marcus's result on multiple  mixing of horocyclic flows \cite{Ma}. Such  flows have  Lebesgue spectrum, so it seemed unsurprising that the flows  possess multile mixing (this notion has been introduced in \cite{Ro} by Rokhlin).  In 1987, I generalized Marcus's theorem as follows: \it if all automorphisms $T_t, \ t>0$  included in an ergodic flow $\{T_t\}$ are isomorphic to each other, then the flow possesses mixing of all ordrs. \rm But this theorem has no interesting applications other than homogeneous flows.  A famous mathematician told me that \it theorem without examples doesn't make much sense.   \rm Years passed, and suddenly the mathematician told me about his theorem. I asked him: Are there any meaningful examples?   \it Theory must be ahead of practice! \rm  That was his clear answer. 

Regarding the actions of Lie groups, I'll point out a difficult problem: \it If the center of action of a Heisenberg group mixes, will it mix multiply? \rm

As for joinings, I rediscovered them in 1986
in connection with Rokhlin's problem on multiple mixing and Furstenberg's theorem on the convergence
$$\sum_{i=1}^n \mu(A_0\cap T^i A_1\cap T^{2i}\cap\dots\cap T^{pi}A_p)/n \ \to
\mu(A_0) \mu(A_1)\dots \mu(A_p)$$
for   weakly mixing  $T$.
Here a  joining $\nu$  appeared   with  the additional invariance  to
$$(T\times T\times\dots\times T)\,\nu\,=\, \nu, $$  it had the following  invariance
$$(Id\times T\times T^2\times\dots\times T^p)\,\nu\, =\,\nu,$$
from which, by induction, it followed that $\nu=\mu\times\mu\times\dots\times\mu$.
I used also a similar additional invariance in the proof of the aforementioned theorem on flows $\{T_t\}$ consisting of isomorphic automorphisms of $T_t$ for $t>0$. So was  my situation befor I looked  at Kalikow's result that had certain philosophical significance. 

Among some experts in ergodic theory, the idea was circulating that slow mixing might not yield multiple mixing. Kochergin flows were mentioned as candidates, \cite{Ko}. Automorphisms of rank 1 have cyclic approximation, so
they have poor statistical properties. All ergodic transformations with positive entropy,  in contrast, require an asymptotically infinite number of cycles in their cyclic representation. And in the case of rank 1, we have only one essentional appoximating cycle.
It is natural to consider  multiple mixing as very good mixing, embedded in a linear scale of mixing properties.
In fact, there is no linear order in this scale! Kalikow's theorem hints that very slow 2-mixing is compatible with 3-mixing. And after the publication of Host's remarkable paper \cite{Ho}, it became clear that this story is completely non-trivial: the worse  2-mixing, the better the 3-mixing it implies. Ledrappier's counterexamples \cite{Le} for the actions of the group $\Z^2$ partially confirm this thesis. They have excellent
2-mixing (Lebesgue spectrum), and the multiple mixing is maximally bad.
From this maximally bad 3-mixing, Tikhonov deduced the triviality of the centralizer for triangular Ledrappier actions \cite{Ti}.

\section{Joinings} 

We consider probability Lebesgue space $(X,\mu)$.
An automorphism (a measure-preserving invertible transformation)
$T:X\to X$ is said to be  {\it k-mixing}, if for $n_1, \dots  n_{k-1}\to+\infty$
$$ \mu(A_0\cap T^{n_1}A_1\cap T^{n_1+ n_2}A_2\cap  \dots \cap T^{n_1+\dots + n_{k-1}}A_{k-1} ) \to  \mu(A_0)\mu(A_1)\dots \mu(A_{k-1}).$$

An automorphism $T:X\to X$  is said
to be of  {\it rank one}, if  there is a sequence $\xi_j$ of measurable partitions
of $X$ in the form
$$\xi_j= \{ E_j,\    TE_j, \  T^2 E_j,\ \dots,   T^{h_j}E_j,  \tilde{E}_j\}.$$
 such that $\xi_j$ converges to the partition onto points.

A {\it self-joining} (of order 2) is defined to be a $T\times T$-invariant
measure $\nu$ on $X\times X$ with the marginals  equal to  $\mu$:
$$\nu(A\times X)=\nu(X\times A)=\mu(A).$$
The joining $\nu$ is called ergodic, if the dynamical system
$(T\times T, X\times X, \nu)$ is ergodic.
We say that $T$ is {\it mixing}, if
$$\Delta^i\to \mu\times\mu,\ \ i\to\infty,$$
i.e.  for all measurable $A,B$
$$\Delta^i(A\times B)=
\mu(A\cap T^iB) \to \mu\times\mu(A\times B) = \mu(A)\mu(B).$$
\vspace{3mm}

We recall that a joining of order $3$ is a $T\times T\times  T$-invariant measure on
$X^3$ such  that all its  margimal projections onto the edges of $X^3$ equal to  $\mu$.

The order $m$ self-joining $\nu$   is called off-diagonal if
it can be represented in the form
$$\nu(A_1\times\dots \times A_m)=
\mu(T^{k_1}A_1\cap \dots \cap T^{k_n}A_n).$$

The property $MSJ(n)$ means that each $n$-fold ergodic self-joining  is
a product of off-diagonals measures.
One can  see that  $MSJ(n)$ is equivalent to  $MSJ(2)$
with the property:
the product measure  $\mu_{(1)}\times \dots \times\mu_{(n)}$
is a unique pairwise independent self-joinings  of order $n$.
Pairwise independent joinings are naturally arising by  the studying
of the   multiple mixing property.

For rank one systems all ergodic joinings have some statistical origin.

\vspace{3mm}
\bf Theorem 1. \it Let  $T$ be a rank-one transformation and $\nu$
be an ergodic self-joining of order 3.
Then there is a sequence $(k1_j,k2_j)$ such that
$(Id\times T^{k1_j}\times T^{k2_j})\Delta \to \frac{1}{3}\nu +\frac{2}{3}\nu'$
for some self-joining $\nu'$.
\medskip

\bf J. King. \it  If  $T$ is a 3-fold mixing and rank-one transformation,
then  $T\in MSJ(3)$.  \rm
\medskip
\medskip

To prove this the following lemma is used (see \cite{92}).

\medskip
\medskip

\bf Choice Lemma. \it Let  $\nu$ be an ergodic measure of \ $T\times T$.
Let a sequence of measures $\nu^k_j$  satisfy the conditions:
for all $A,B$
$$\left|\nu^k_j(TA\times TB) - \nu^k_j(A\times B)\right| < d(j)\to 0$$
and  $$\sum_k a^k_j  \nu^k_j \to \nu \ \   (\sum_k a^k_j =1).$$
Then there is  a sequence $k_j$ such that
$  \nu^{k_j}_j \to \nu.$

In other words,
if the average of almost invariant measures
is close to an ergodic measure, then there exists a weighted majority of measures 
in which every measure is close to this ergodic measure. \rm

\rm

\vspace{3mm}
Proof of lemma.  Given sets $A, B$ and  $\eps >0$
we consider the sets $K_j \subset [-\delta h_j,\delta h_j]$ of all integers $k$
such that
$$ \nu(A\times B) -  \nu^k_j  (TA\times TB) >\eps$$
( or
$ \nu^k_j  (A\times B)-\nu(TA\times TB)  >\eps \ )$).
Suppose that  the (sub)sequence  $K_j$
satisfies the condition
$$\sum_{k\in K_j} a^k_j \geq a>0.$$
Let $\lambda$ be a limit point for the sequence  of the measures
$(\sum_k a^k_j)^{-1}\sum_k a^k_j  \nu^k_j $. The measure $\lambda$  is invariant,
and $\lambda\neq \nu$.
However we  represent ergodic measure $\nu$ as
$$\nu= a\lambda +(1-a)\lambda',$$
for   some invariant  measure $\lambda'$.  From ergodicity of $\nu$ we get
$\lambda=\nu$. The contradiction shows that
$$\sum_{k\in K_j} a^k_j \to 0.$$
Thus, for most of the measures $\nu^k_j$ we have
$$| \nu^k_j  (A\times B) - \nu(A\times B) |< \eps.$$
For a family $\{A_1, A_2,\dots\}$, which is dense in the algebra of
all $\mu$-measurable sets,
using diagonal method we find a sequence  $k_j$ such that
$$| \nu^{k_j}_j (A_m\times A_n) - \nu(A_m\times A_n) | \to 0.$$
This implies
$$| \nu^{k_j}_j (A\times B) - \nu(A\times B) | \to 0$$
for all $A,B$.

\section{ Lemmas for mixing rank one transformations}         \rm

\bf  Blum-Hanson's lemma. \it Let $T$ be mixing.
Suppose the sequence $\lbrace a^{z}_j \rbrace , z\,\in\,Z ,j\,\in\,\bf N$,
satisfies the conditions:
$$
\sum_z a^{z}_j =1, a \geq 0,$$
 $$ max_z \lbrace a^{z}_j \rbrace
\,\to \,0,\,j\,\to \,\infty.$$
Then we have
$$
  \left\| \sum_z a^{z}_j T^zf - \Theta \! f\, \right\|_2 \,\to \,0.
$$
 \rm
\vspace{3mm}

Proof.  Let $\Theta  f=0$.  Put $P_j = \sum_z a^{z}_j T^z.$
Let us show that $  \| P_jf\|_2 \,\to \,0.$
One has
$$P_j^\ast P_j = \sum_w b^{w}_j T^w ,$$
where  the sequence $\lbrace b^{w}_j \rbrace$ satisfies
$$   b_{w}^j \leq \sum_z a^{w-z}_j a^{z}_j \leq max_z {a^z_j } \to \,0.
$$
Since $T$ is mixing, one has   $\sum_w b^{w}_j T^w f \to 0$ (weakly).
Thus,
$$\| P_j f\|^2=(P_j^\ast P_j f, f)\to 0.$$

If $\Theta  f\neq 0$, we get
$$\|P_j(f-\Theta f)\|^2 \to 0, \ \ \ \|P_jf-\Theta f\|^2\to 0.$$

\vspace{3mm}

\bf Lemma on small returns. \it Let $T$ be a mixing rank-one transformation  and 
$a^{z}_j: =\ \mu(T^zE_j | E_j)$.
Then $\lim_j max_{z>0} \{a^z_j\} = 0$.  \rm
\vspace{3mm}

Proof. We have  $ max_{z>0} \{a^z_j\}= max_{z>h_j} \{a^z_j\}$.
Suppose $\lim_j max \{a^z_j :z>h_j\} = a>0$,
 $\mu(T^{z_j}E_j | E_j)\to a,$ hence,
$$\mu(T^{z_j}T^kE_j | T^kE_j)\to a, \ \ (0\leq k\leq h_j).$$
We can approximate the measurable set $A$  by $\xi_j$-measurable
sets $A_j$ ( $A_j$ is a union of some $T^kE_j$). We have
$$\limsup_j \mu(T^{z_j}A_j | A_j)\geq  a,$$
this implies
$$\limsup_j \mu(T^{z_j}A | A)\geq  a.$$
By $\mu(A)<a$ we obtain the contradiction with the property of mixing:
$\mu(T^{z_j}A | A)\to \mu(A)$.  Thus,  $$\lim_j max \{a^z_j :z>0\} = 0.$$

\section{   Spraying pairwise independent joinings}

\bf The main  theorem. \it Let $T$ be a mixing rank-one transformation.
Let $\nu$ be a pairwise independent self-joining of order $ 3$ (the projections of $\nu$ onto the two-dimensional faces of the Cartesian cube coincide with  $\mu\times \mu$.
Then $\nu=\mu^3:=\mu\times \mu\times \mu$.  \rm

\vspace{3mm}
For $\eps >0$ we define $\eps$-light blocks:
$$
  D_{\eps, j}=\{z\in [0,h_j]^3 \ :\  z=(z1,z2,z3),
\bar{E}^z_j= T^{z1}E_j\times  T^{z2}E_j\times  T^{z3}E_j
\nu(\bar{E}^z_j)<\eps \mu(E_j)^2\}.
$$
Now we calculate the total mass $Di(\nu)$ of "infinite light blocks":
$$
  Di(\nu)=\lim_{\eps\to 0}
\left(\limsup_{j\to \infty} \sum_{z\in D_{\eps,j}}\nu(\bar{E}^z_j)\right).
$$
\vspace{3mm}

\bf Heavy blocks generate  light blocks.  \it
If $\nu$ is a pairwise independent self-joining,
then $Di(\nu)>0$.  \rm
\vspace{3mm}

We can see that  small returne  lemma show that   the heavy blocks
under the action of  some powers of $T\times T\times T$  generate
many light blocks. So   $Di(\nu)>0$.

\begin{center}
\includegraphics{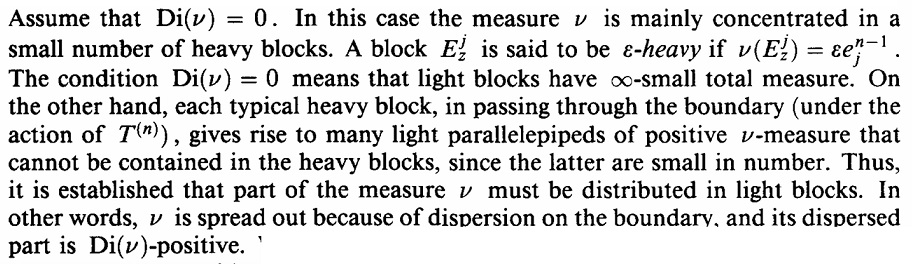}
\end{center}

\bf Dust destroys non-triviality of  self-joinings. \it
Let  $\nu$ be  a pairwise independent ergodic self-joining,
 and
 $Di(\nu)>0$, then $\nu =\mu\otimes \mu\otimes\mu$.
\rm
\medskip

Proof.
Let's define  columns
in the following way:
$$
C_j^{w,s}=\uu_{i=0}^{\delta h_j}
T^{w+i}E_j\times T^{i}E_j\times T^{s+i}E_j.
$$
Given small $\delta >0$ we   find  a sequence of
sets  $F_j$ such that
$$\nu(\ | F_j)\to \nu$$
and the sets $F_j$ ( fibers) have the form
$$
F_j =:\uu_{h\in S_j} (Id\times Id\times T^h)C_j,$$
for some sequences
 $$S_j\subset \{0,1, \dots (1-\delta) h_j\},
\ \ s_j, w_j\in \{0,1, \dots (1-\delta) h_j\},
\ \ C_j= C_j^{w_j,s_j}.$$
Since the sets as $F_j$  are almost invariant
 with respect to $T\times T\times T$   and
 $\nu$ is ergodic,  from Choice Lemma we  get  some sequence $F_j$ such that
$\nu(\ | F_j)\to \nu$.

Moreover  we have  $Di(\nu)>0$, so  one can find  a sequence of fibers
$F_j$
with  such  sets $S_j$ that numerate only light columns.
The sets $S_j$ will  satisfy the condition:
$$\max_{s\in S_j}\{a_s^j\} \to 0, \  j\to\infty, \ \ \
 \sum_{s\in S_j}a_s^j =1,
$$
where
$$a_s^j=\nu((Id\times Id\times T^s)C_j   \,|\, F_j).$$

Let   $Y_j$  be the projection of $C_j$   into the factor $X_{(2)}$ of
the product
 $X_{(1)}\times X_{(2)}\times X_{(3)}$.
The measure
$\nu(\ | F_j)$  is close to $\lambda_j$ (a part of off-diagonal measure
in $X\times X\times X$),  which is defined by the equality
$$ \lambda_j(A\times B \times C)\ =:
\ \frac{1}{a\mu(Y_{j})}
\sum_s a_s^j \int_X \chi_{Y_{j}}\, \chi_{T^{w_j}A}
\,\chi_B\, \chi_{T^s C}\, d\mu.
$$
From Blum-Hanson's theorem we get
$$\sum_{s\in S_j} a_s^j \chi_{T^s C} \ \to_{L_2}  \ Const \
\equiv \ a\mu(C) \   (j\to\infty).
$$
Using  pairwise independence of the self-joining $\nu$ we obtain
$$\nu(A\times B \times C) =
\lim_{j\to\infty} \nu(A\times B \times C\, |\, F_j) =
 \lim_{j\to\infty}  \lambda_j(A\times B \times C)
$$
$$
=\lim_{j\to\infty} \ \frac{1}{a\mu(Y_j)}
\int_X \chi_{Y_j}\, \chi_{T^{w_j}A}   \,\chi_B\,
\left( \sum_{s\in S_j} a_s^j   \chi_{T^s C}\right)\, d\mu
$$
$$ =
\mu(C)\,\nu(A\times B \times X) =
\mu(A)\mu(B)\mu(C).
$$
\Large

Kalikow's theorem left open the possibility of generalizations not only to high multiplicities, but also to actions of groups $\Z^n$, $\R^n$, $n>1$ (here the reader might want to think deeply about what rank 1 means for actions of these groups). One can also attempt to prove multiple mixing by weakening the condition "be an action of rank 1" {}, requiring instead that the rank be finite or the local rank be positive.
Surprisingly, much of this program has been implemented; see \cite{92},\cite{93},\cite{2000}.

King proved that mixing rank 1 automorphisms have the $MSJ(2)$ property,
and due to  Kalikow's theorem, they have the $MSJ(3)$ property.

Recall that Ornstein proposed a stochastic ensemble of rank-1 transformations that, with probability one, have the  trivial centralizer.
Remarkably, the approximation and mixing properties here imply an algebraic property. Examples of the surprising connection between algebra and statistics exist beyond rank 1 as well: Austin's result \cite{Au} on the weak Pinsker property implies that a K-automorphism (an automorphism with completely positive entropy) is a direct product of countably many automorphisms, so the centralizer of every K-automorphism is incredibly huge.

Rudolph's work \cite{Ru} made a big impression by providing an example of a rank-1 mixing transformation with minimal self-joinings of all ordes ($MSJ$), demonstrating spectacular applications of the $MSJ$ property. Invent a mixing automorphism that has a root of degree 2026 but no root of degree 2027. Rudolph easily constructed such toys using MSJ-bricks. The bricks themselves have the trivial structures of invariant sigma algebras (factors) and the
trivial centralizer. A distinctive feature of  Rudolph's automorphism $R$ is that
a set of distinct factors isomorphic to $R$, for any ergodic system containing them, is a globally independent set of facts. $R$-factors are either completely glued together or independent. All this made it possible to completely describe the structures of factors and centralizers of direct products of transformations with  $MSJ$-property.

The topic of minimal self-joinings has been stagnant for a long time, but a scientific feat was recently achieved: the puzzling technique of Chaika and Roberson \cite{CR} allowed us to solve Thouvenot's problem: \it

Mild mixing (no rigid factors) with rank 1 together do not imply minimal self-joinings. \rm

King \cite{Ki} once unexpectedly showed that $MSJ(4)=MSJ:=MSJ(\infty)$. This is one of my favorite assertions; it is a clever consequence of a completely trivial but very useful fact about bounded operators on Hilbert spaces:
$$A^\ast A=0 \ \Rightarrow \ A=0.
$$
For commutative actions
$MSJ(3)=MSJ$, see \cite{GHR}, \cite{R8} (for non-commutative ones, there are counterexamples!).
For automorphisms with singular spectrum and finite-rank $\Z^n$-actions, $MSJ(2)=MSJ$. For flows without any restrictions, we have $MSJ(2)=MSJ$ (in   more general situations  \cite{RT} for  applications see  \cite{FK}).

$\Z^n$-actions with the 3-fold mixing and the property of $MSJ(2)$ have
$MSJ=MSJ(3)$ and multiple mixing of all orders and mixing of all orders. Over the past 30 years, there has been no progress coserning  the problem $MSJ(2)=MSJ(3)$? If there is a counterexample with simple  spectrum, then  this example  solves both the Banach simple Lesbegue spectrum  problem  and Rokhlin multiple mixing problem.

\vspace{5mm}
\bf Spectral associations with  Kalikow's $\bf (T,T^{-1})$-automorphism. \rm   Changing the plot, I can offer a problem related to another well-known work by Kalikow, \cite{Kalikow}.

Given a sequence of numbers
$$s_{2n-1}=0,  \ \ \ s_{-2n}=  s_{2n}, \ \ 
s_{2n}=2^{-2n}\, C_{2n}^{n}, \ \ n\geq 0,$$  show that
$s_n$ are the Fourier coefficients of some continuous measure $\sigma$ on the circle,
$s_n=\widehat{\sigma(n)}$.   Then prove the absolute continuity of this measure, using and not using  Kolmogorov's theorem on spectra of K-automorphisms.
  
\large


\begin{thebibliography}{99}
 

\bibitem{Ka} S.A. Kalikow, Two fold mixing implies threefold mixing for rank one transformations, Ergod. Th. Dynam. Sys., 4 (1984), 237-259

\bibitem{Ro} V.A. Rokhlin, On endomorphisms of compact commutative groups, Izv. Akad. Nauk SSSR Ser. Mat., 13:4 (1949), 329-340

\bibitem{Ma} B. Marcus, The horocycle flow is mixing of all degrees, Inv. Math., 46 (1978), 201-209 

\bibitem{Ko}
 A.V. Kochergin, On mixing in special flows over a shifting of segments and in smooth flows on surfaces, Math. USSR-Sb., 25:3 (1975), 441-  469

\bibitem{Th}  J.-P. Thouvenot, Some properties and applications of joinings in ergodic theory, Ergodic Theory and Its Connections with Harmonic Analysis, Proceedings of the Alexandria 1993 Conference, LMS Lecture Note Series, 205, eds. K. E. Petersen and I. A. Salama, 1995, 207-235  



\bibitem{Au}	 T. Austin, Measure concentration and the weak Pinsker property, Publ. Math. IHES, 128:1 (2018), 1-119 

\bibitem{Le}F. Ledrappier, Un champ marcovien peut \^ etre d'entropie nulle et m\' elangeant, C. R. Acad. Sci. Paris S\' er. A-B, 287:7 (1978), A561-A563

\bibitem{Ti}S.V. Tikhonov, On the Absence of Multiple Mixing and on the Centralizer of Measure-Preserving Actions, Math. Notes, 97:4 (2015), 652-656 

\bibitem{92}
V.V. Ryzhikov, Mixing, rank, and minimal self-joining of actions with an invariant measure, Russian Acad. Sci. Sb. Math., 75:2 (1993), 
405-427  

\bibitem{93}V.V. Ryzhikov, Joinings and multiple mixing of finite rank actions. Funct. Anal. Appl. 27 (1993), No. 2, 128-140

\bibitem{2000}
V.V. Ryzhikov, Rokhlin's multiple mixing problem in the class of positive local rank actions, Funct. Anal. Appl., 34:1 (2000), 73-75  


\bibitem{KT} 	J.L. King, J.-P. Thouvenot, A canonical structure theorem for finite joining-rank maps, J. Analyse. Math., 56 (1991), 211-230 

\bibitem{Ho} B. Host, Mixing of all orders and pairwise independent joinings of systems with singular spectrum, Israel J. Math., 76:3 (1991), 289-298  

\bibitem{Ki}	J.L. King, Ergodic properties where order 4 implies infinite order, Israel J. Math., 80:1-2 (1992), 65-86  

\bibitem{GHR} E. Glasner, B. Host, D. Rudolph, Simple systems and their higher order self-joinings, Israel J. Math., 78 (1992), 131-142 


\bibitem{R8} V.V. Ryzhikov, Self-Joinings of Commutative Actions with Invariant Measure, Math. Notes, 83:5 (2008), 723-726


\bibitem{97}V.V. Ryzhikov, Around simple dynamical systems. Induced joinings and multiple mixing. J. Dyn. Control Syst. 3 (1997), No. 1, 111-127



\bibitem{Ru}D. Rudolph, An example of measure-preserving map with minimal self-joinings, and applications, J. d'Analyse Math., 35 (1979), 97-122 

\bibitem{JR}A. del Junco, D. Rudolph, On ergodic action whose self-joinings are graphs, Ergodic Theory Dynam. Systems, 7 (1987), 531-557


\bibitem{CR}Jon Chaika, Donald Robertson, A rank one mild mixing system without minimal self joinings, Annales Henri Lebesgue, 9 (2026), 97


\bibitem{RT}  V.V. Ryzhikov, J.-P. Thouvenot, Disjointness, divisibility, and quasi-simplicity of measure-preserving actions, Funct. Anal. Appl., 40:3 (2006), 237-240 


\bibitem{FK}  B. Fayad, A. Kanigowski, Multiple mixing for a class of conservative surface flows, Invent. Math., 203:2 (2016), 555-614 

\bibitem{Kalikow}  S.A. Kalikow. $T,T^{-1}$-transformation is not loosely Bernoulli. Ann. of  Math. (2), 115(2):393-409, 1982
\end{thebibliography}
\end{document}